\documentclass{amsart}[11pt]
\usepackage{myAMSpackages}
\usepackage{Universal}
\usepackage{MyAMSEnvironments}
\usepackage{graphicx,subfigure,epsfig}

\newcommand{\wh}[1]{\widehat{#1}}

\newcommand{\tH}{\ensuremath{\text{\tiny $H$\normalsize }}\xspace}

\newcommand{\J}{\mathcal{J}}

\newcommand{\bp}{\begin{pmatrix}}
\newcommand{\ep}{\end{pmatrix}}

\newcommand{\bphi}{\boldsymbol \varphi}
\newcommand{\E}{{\boldsymbol E}}
\newcommand{\F}{ {\bf F}}
\newcommand{\Z}{ {\boldsymbol Z}}
\newcommand{\ee}{ {\boldsymbol e}}

\begin{document}
\title{Intrinsic complements of equiregular sub-Riemannian manifolds}
\author{Robert K. Hladky}

\begin{abstract}
Under a nondegeneracy condition, we show that an equiregular sub-Riemannian manifold of step size $r$ admits a canonical, $V$-rigid complement defined from the sub-Riemannian data that is preserved the by action of sub-Riemannian isometries. We explore how the existence of such a complement relates to results from the literature and study the step size $2$ case in more detail.
\end{abstract}

\maketitle

\section{Introduction}

In recent years, there has been activity in applying tensorial analysis on sub-Riemannian manifolds. Topics studied include minimal and constant mean curvature spaces (\cite{CHMY}, \cite{HP}, \cite{HP2} ,\cite{HP4},  \cite{Montefalcone}), natural sub-elliptic operators and heat kernels methods (\cite{BaudoinGarofalo} , \cite{BaudoinGarofalo2} ,\cite{BBG} ,\cite{BBGM}, \cite{Hladky5}), and isometries of sub-Riemannian manifold (\cite{Hladky6}), amongst others.  A key feature of many of these tensorial approaches has been the choice of a complementary bundle and the definition of some type of connection dependent on this bundle. Indeed, in \cite{Hladky4}, the author shows, that given a choice of complement, a canonical notion of covariant derivative can be defined for sections of the horizontal bundle.

In this short paper, we show that, for equiregular sub-Riemannian manifolds satisfying a nondegeneracy condition, there is an intrinsic choice of complementary bundle together with an extension of the sub-Riemannian metric to a full  Riemannian metric. The complement and metric extension are  also both invariant under smooth infinitesimal sub-Riemannian isometries. The canonical connection associated to the complement and metric extension defined in \cite{Hladky4} has the $V$-rigid property (defined in \rfS{obs} or \cite{Hladky4}) which, as we shall discuss, greatly simplifies tensorial analysis in many applications.

In \rfS{def} we develop the required notation and terminology, briefly review known intrinsic objects on equiregular sub-Riemannian manifolds and introduce the bracket-structure $\J$-maps and the semi-$\J$-nondegenerate condition. In \rfS{proof}, we prove our main result on the existence of an intrinsic natural complement. 

\begin{thms}\label{TS:A}
For any semi-$\J$-nondegenerate, equiregular sub-Riemannian manifold there exists an algorithmically defined,  natural , $V$-rigid complement.
\end{thms}
The complement produced by this algorithm is named the {\em minimal rigid complement}.

The precise notion of complement and vertical rigidity are introduced at the end of \rfS{def}. The adjective natural here means that the object is invariant under smooth infinitesimal sub-Riemannian isometries. In \rfS{obs} we remark on the relevance of this result to several topics that are prevalent in the literature. In particular, we establish

\begin{thms}If $M$ is a a semi-$\J$-nondegenerate, equiregular sub-Riemannian manifold, the group of smooth, infinitesimal sub-Riemannian isometries of $M$ is a Lie group with dimension bounded above by $\frac{\dim H (\dim H-3)}{2} + \dim M    $.
\end{thms}
 Finally in \rfS{two}, we discuss the step $2$ case where the nondegeneracy condition is least restrictive.

\section{Terminology and definitions}\label{S:def}

In this section, we introduce the basic definitions and terminology used throughout this paper.

\bgD{subR}
A sub-Riemannian manifold $(M,H,g_\tH)$ consists of a smooth manifold $M$, a smooth bracket-generating subbundle $H$ of the tangent bundle and a smoothly varying, positive definite inner product $g_\tH$ on each $H_x$ with $x \in M$. 
\enD

We shall focus on equiregular sub-Riemannian manifolds which have a strong additional assumption on the Lie brackets of vector fields in $H$. To describe this assumption, let $H_1 = H$ and define a nested sequence of distributions
\bgE{Hseq}  H_1 \subseteq H_2 \subseteq H_3 \subseteq \dots \enE
by 
\[ (H_{i+1})_x = (H_i)_x + \{ [A,B]_x \colon A \in C^\infty( H_i), B \in C^\infty(H_1) \}  .\]
The bracket-generating condition can then be interpreted as for all $x\in M$, there is an integer $r(x)$ such that $(H_{r(x)})_x  =T_x M$. We say that $M$ has step size $r$ if $H_r =TM$ and $H_{r-1} \ne TM$.
\bgD{equiregular}
A sub-Riemannian manifold $M$ is equiregular if it has step size $r$ for some $r \geq 1$  and for each $i=1,\dots, r$ the distribution $H_i$ is of constant rank.
\enD
A trivial consequence of equiregularity is that the inclusions $H_i \subset H_{i+1}$  in \rfE{Hseq} are strict for $i<r-1$ and equalities for $i \geq r$.  For the remainder of this paper, we shall always assume that $M$ is equiregular with step size $r$.

Next we introduce an intrinsic way of encoding the differences between these successive subbundles. Set $H_0 = \{0\}$ and then inductively  $\wh{H}_m = H_m / H_{m-1}$. For convenience of notation, we set $d_m = \dim \wh{H}_m$ and note $\dim M = d_1 + \dots + d_r$. For $A \in (H_m)_x$, let $[A]_m$ denote the equivalence class in $\wh{H}_m$. Now we introduce the annihilator subbundles of $T^*M$ by
\[  (H^o_0)_x = T^*_xM , \qquad (H^o_i)_x = \{ \phi \in T_x^*M \colon \phi_{| (H_i)_x } = 0 \}.\]
Thus $H_r^o  \subset H^o_{r-1} \subset \dots \subset H^o_2 \subset H^o_1 \subset  (H^o_0)= T^*M$ and we can define natural dual spaces to $\wh{H}_1,\dots,\wh{H}_r$ by setting $\wh{V}^m = H^o_{m-1} / H^o_{m}$ for $m=1,\dots,r$.

The geometry of a sub-Riemannian manifold is closely connected to the bracket-structure on $H$.  This structure can be encoded in a family of pointwise bilinear mappings
\[ B^{k,m} \colon \wh{H}_m \times \wh{H}_k \to \wh{H}_{k+m} \]
defined as follows: if $a \in (\wh{H}_k)_x$,  $b \in (\wh{H}_m)_x$, let $A,B$ be any sections of $H_k$, $H_m$ respectively that represent $a,b$ at $x$, then 
\bgE{B} B^{k,m}(a,b) =\big[  -[A,B]_{|x} \big]_{k+m}\enE

\bgL{welldef}
The bilinear forms $B^{k,m}$ are well-defined.
\enL

\pf
Given sections $A,B$ of $H_k$, $H_m$ representing $a$, $b$ respectively. Any other such sections representing $a$ and $b$ will then locally take the form 
\begin{align*}
\tilde{A} &= A + f^i Z_i + C \\
\tilde{B} &= B + g^j W_j +D
\end{align*}
where $\{Z_i\}$, $\{W_j\}$  are sections of $H_k$, $H_m$ that represent orthonormal frames for $\wh{H}_k$, $\wh{H}_m$ respectively, $f^i$ and $g^j$ are collections of smooth functions that vanish at $x$ and $C$, $D$ are sections of $H_{k-1}$, $H_{m-1}$ respectively.

Then, using the fact that each $f^i,g^j$  vanishes at $x$ and $[\tilde{A},D]_{|x} , [ C,\tilde{B}]_{|x} \in H_{m+k-1}$, we must have
\begin{align*}
[ \tilde{A}, \tilde{B} ]_{|x} &= [A,B ]_{|x} + [f^iZ_i,\tilde{B}]_{|x} + [\tilde{A},D]_{|x} + [C,\tilde{B}]_{|x} + [\tilde{A},g^jW_j ]_{|x}  \\
&= [A,B ]_{|x}   -  \left( (\tilde{B} f^i) Z_i \right)_{|x}  + [\tilde{A},D]_{|x}+[C,\tilde{B}]_{|x} + \left( (\tilde{A} g^j) W_j  \right)_{|x} \\
&= [A,B]_{|x}   \qquad \text{mod}  (H_{k+m})_{|x} .
\end{align*}

\epf

Our purpose is to construct a complementary bundle to $H$ for an equiregular sub-Riemannian manifold that is intrinsic in the sense that it depends only on the sub-Riemannian data. On sub-Riemannian manifolds with plentiful symmetries, one consequence is that the complement must be invariant under appropriately defined sub-Riemannian isometries. Following the language of \cite{Hladky6}, we make the following definition.

\bgD{weakH}
If $(M,H,g_\tH)$ is a sub-Riemannian manifold, a weak $H$-isometry of $M$ is a diffeomorphism $F \colon M \to M$ such that $F_* H = H$ and $F^* g_\tH = g_\tH$.
\enD

We shall call an object {\em natural} if it is preserved in an appropriate sense by all weak $H$-isometries of $M$ in an appropriate sense. As immediate examples, it is trivial to show that the distributions $H_m$ and $H_m^o$ are natural as $F_* H_m = H_m$ and $F^* H_m^o = H_m^o$ for all diffeormorphisms preserving $H$. For weak $H$-isometries it is then  easy to verify that the bilinear mappings $B^{k,m}$ are natural in the sense that $F_* B^{k,m}(a,b) = B^{k,m}(F_*a,F_*b)$.

We now show that the inner product can be extended in an intrinsic fashion to the quotient bundles $\wh{H}_m$. Define $ \pi_1 \colon  \wedge^2 H_1 \to \wh{H}_2 $
by linear extension of 
\[ \pi_1 ( X \wedge Y) =  B^{1,1}(X,Y). \]
For $j>1$,  we define linear maps $\bt{\pi}_j \colon H_1 \otimes \wh{H}_j \to \wh{H}_{j+1}$ by linear extension of $B^{1,j}$.  Then for all $j>0$, we can inductively identify  $\wh{H}_{j+1}$ with $\ker{\bt{\pi}_j}^\perp$. We then use this to successively  define an inner product on each $\wh{H}_m$.  By duality, we then also have inner products on each $\wh{V}^m$.  It is easy to verify that each of these inner products is natural in the obvious sense.

The key difficulty with using the bundles $\wh{H}_m$ is that sections are not genuine vector fields but sections of $H_m$ defined only up to the addition of a section of $H_{m-1}$. We shall need language to describe particular representatives in each equivalence class.

\bgD{rep}
For $m>k$, an $(m,k)$-frame at $x$  is a row-vector $\E_x =(E_1,\dots,E_{d_m})$  with entries in $\left(H_m / H_k \right)_x$ such that  the equivalence class $[\E_x]_m$ forms an orthonormal frame for $( \wh{H}_m)_x$. Thus an $(m,m-1)$-frame field consists of sections of $\wh{H}_m$, but an $(m,0)$-frame field consists of genuine vector fields.

An $m$-coframe at $x$ is a column-vector  $\bphi=(\phi^1,\dots, \phi^{d_m})$ of elements of  $(H^o_{m-1})_x$ such the equivalence class  $[\bphi]^m$ is an orthonormal frame for $\wh{V}_x^m$.

A pair $(\E_x,\bphi_x)$ consisting of an $(m,k)$-frame and $m$-coframe at $x$ is dual if $\bphi(\E) = I_{d_m \times d_m} $.
\enD

Let $\Phi^m_x $ denote the collection of $m$-coframes,  $\bphi_x$. For $\bphi_x \in \Phi_x^m$ and $k<m$, let $\Omega^{m}_k(\bphi_x)$ denote the collection of $(m,k)$-frames at $x$ dual to $\bphi_x$. Thus for a fixed element $\bphi_x \in \Phi_x^m$, the set $\Omega^m_{m-1}(\bphi_x)$ consists of a single element. Furthermore,  each element of  $\Omega^m_k(\bphi_x)$ can be viewed as an affine space acted on by  $ \left(H_{k-1} \right)_x^{d_m}$. We shall typically work with local sections of these spaces. 

It is also easy to check that under an orthonormal change of frame $f  \bphi  $, 
\[ \Omega^m_k( f \bphi )  = \Omega^m_k(\bphi) f^\top \]
where $ [\E]_k f^\top = [\E f^\top]_k$. Thus we can also define pointwise
\[ \Omega^m_k = \bigcup\limits_{ \bphi \in \Phi^m} \Omega^m_k / \sim \]
where $\sim$ is equivalence under a orthonormal change of frame, i.e.  $[\E] \sim [\E f^\top]$ for any section $f$ of $O(d_m)$.

The following observations are trivial.
\bgL{affine}  
If $\E$ is an $(m,k)$-frame at $x$, then $\E$ is naturally isomorphic to a set $(m,k-1)$-frames at $x$ that form an affine space over $(\wh{H}_{k-1}^{d_m})_x$.
\enL

\bgL{OmNat}
The bundles $\Omega^m_k$ are natural.
\enL

\pf This follows from the trivial observation that for a weak $H$-isometry $F$,
\[ F_* \Omega^k_j(F^* \bphi) = \Omega^k_j(\bphi).\]

\epf

We can now introduce our nondegeneracy condition. Using the intrinsic inner products on each $\wh{H}_i$, we obtain well-defined, natural, linear bundle maps
\[ \J^{m,k}  \colon \wh{V}^{m+k} \to \text{Hom}(\wh{H}_m, \wh{H}_k) \]
by
\bgE{J} \aip{\J^{m,k}([\phi])a}{b}{}  :=  [\phi] ( B^{m,k} (a,b) ) = d\phi(a,b).\enE

\bgD{sndg}
An equiregular sub-Riemannian manifold $M$  is  semi-$\J$- nondegenerate at $x$ if for all $2\leq k \leq r$ and   $k$-coframes $\bphi_x \in \Phi^{k} (\bphi_x)$, the map \[ \widehat{\J}^{k}_{\bphi_x} \colon  \left( \wh{H}_{k-1} \right)_x^{d_{k}} \to M_{d_{k} \times d_{k}}(H_x)\]  defined by 
\bgE{wJ}\widehat{\J}^{k}_{\bphi_x} (\Z) =\J^{k-1,1} {[ \bphi_x ]}(\Z) +  \left[\J^{k-1,1}{[\bphi_x]}(\Z )\right]^\top\enE is injective. We say $M$ is semi-$\J$-nondegenerate if this holds at every $x\in M$.
\enD

Here for a vector space $Y$, the notation $Y^k$ denotes the vector space of $k$-row vectors with entries in $Y$ and $M_{m \times n}(Y)$ is the space of $m \times n$ matrices with entries in $Y$. In the step $2$ case discussed in \rfS{two}, there is a stronger nondegeneracy condition called $\J$-nondegenerate which motivates the modifier "semi" used here.

It should be remarked that under an orthonormal change in $\bphi$,
\[ \widehat{\J}^{m}_{f \bphi}(\Z f^\top) = f \widehat{\J}^{m}_{\bphi} (\Z) f^\top \]
and so the semi-$\J$-nondegenerate condition is independent of the choice of $\bphi_x$.

As the final part of this section, we introduce a precise definition of the notion of a complement.

\bgD{complement}
A partial complement at step $m+1$  is a subbundle such $V_{m+1}$
\[ H_{m+1} = H_m \oplus V_{m+1} \]
for $m=1,\dots, r-1$.
A graded complement $V$ is a direct sum of partial complements $V=V_2 \oplus \dots \oplus V_r$.
 \enD
 
These are the first objects under consideration in this paper that are not defined from intrinsic sub-Riemannian data.  It should then be emphasized here that partial complements typically are not natural.  Indeed our primary focus is to establish the existence of a natural complement under the semi-$\J$-nondegenerate condition.

However, any choice of partial complement $V_j$ is pointwise isomorphic to $\wh{H}_{j}$ and so inherits an inner product. Thus any choice of complement comes equipped with an intrinsic extension of the sub-Riemannian inner product with the property that
\[ TM = H \oplus V_2 \oplus \dots \oplus V_r \]
is an orthogonal decomposition.

 For many applications a further property of complements is often useful.
  \bgD{srigid}
 A complement $V$ is $V$-rigid if whenever $\{U_i\}$ is a local orthonormal frame for $V$ with dual frame $\{\psi^i\}$, we have
 \bgE{vra} \sum\limits_i  d\psi^i(U_i,X)  =0 \enE
 for all horizontal vectors $X \in H$.
 \enD
If $f$ is an orthogonal matrix valued function,  $df f^\top$ is easily seen to be a skew-symmetric matrix. Thus this condition need only be checked for one orthonormal frame near each point.   Some consequences of this definition are reviewed in \rfS{obs}.
 
We end this section by remarking that locally a partial complement $V_m$ is spanned by collection of vector fields $E_1,\dots,E_{d_m}$ such that the equivalence classes $[E_1]_m,\dots,[E_{d_m}]_m$ form a local orthonormal frame field for $\wh{H}_m$.   Any other complement will be spanned by alternative representatives from the same equivalence classes.  A choice of partial complement can be regarded as a choice of representatives for $[E_1]_m,\dots,[E_{d_m}]_m$.



\section{Existence of an intrinsic natural connection}\label{S:proof}

In this section we shall prove Theorem A by introducing measures on the spaces of local  $(m,k)$-frame fields $[\E]$   viewed as affine spaces of $(m,k-1)$-frame fields. These measures shall be defined in terms of the norms of matrix-valued functions. Our algorithm will essentially run by iteratively isolating unique sections of minimal norm until we have a well-defined $(m,0)$-frame for each $m$.

 Before we can begin the proof, we need a few more preliminary definitions and results. First we formally construct the norm we shall be minimizing. Given an inner product space $Y$, we set $M_{j \times k}(Y)$ to be the vector space of $j \times k$ matrices with entries in $Y$. We can then define an inner product on $M_{j \times k}(Y)$ by
\[ \aip{A}{B}{} =  \sum\limits_{l,m}  \aip{(B^\top)_l^m}{A_m^l}  .\]
Informally, this can be expressed as $\aip{A}{B}{} = \text{trace} B^\top A$ where all multiplications in the product are replaced by inner products in $Y$. In particular, $\|A\|^2 = \sum \|A_l^m\|^2$.

Given an $m$-coframe field $\bphi$  and a dual $(m,m-j-1)$-frame field $\E$ with $j\geq 1$, we can define a section $\mathscr{S}(\bphi,\E)$ of $M_{d_m \times d_m}(H)$ by
\bgE{Sdef} \begin{split} \aip{\mathscr{S}(\bphi,\E)}{X}{} &=  d\bphi(\E,X)  + d\bphi(\E,X)^\top.   \end{split}\enE
This is well-defined as $\E$ is uniquely defined up to sections of $H_{m-j-1}$ and  \[[X,H_{m-j-1}] \subset H_{m-1}\] which is annihilated by $\bphi$. 

Given an orthonormal change of frame $f \bphi$ with corresponding change of dual coframe $\E f^\top$, we obtain
\begin{align*}
d(  f \bphi ) ( \E f^\top ,v) &=  f d \bphi(\E,v) f^\top -df(v) f^\top . 
\end{align*}
Now $f df^\top + df f^\top = 0 $ so the later portion is skew-symmetric and hence
\bgE{Sf} \mathscr{S} (f \bphi, \E f^\top) =  f \mathscr{S} (\bphi,\E) f^\top .\enE

Additionally, given an $m$-coframe field $\bphi$ and a $(k,m-1)$-frame field $\E$ with $m < k$ and $\bphi(\E)=0$, we define a section of $M_{d_m \times d_k}(H)$ by
\bgE{Adef}
\aip{\mathscr{A}(\bphi,\E)}{X}{} =  d\bphi(\E,X).
\enE
With these conditions on $\bphi$ and $\E$, it follows that under orthonormal changes to $f_m \bphi$ and $\E f_k^\top$,
\bgE{Af}
\mathscr{A}(f_m \bphi,\E f_k^\top) =f_m \mathscr{A}(\bphi,\E) f_k^\top
\enE
and so
\[ \| \mathscr{A}(f_m \bphi,\E f_k^\top) \| = \| f_m \mathscr{A}(\bphi,\E) f_k^\top\|. \]
Now for $j<m$ an element  $\E \in \Omega^m_{j}(\bphi_x) $ can be viewed as an affine space $[\E]$ with elements in $\Omega^m_{j-1}(\bphi)$ and underlying vector space $(\wh{H}_{j})^{d_m}$. Hence for a $m$-coframe field $\bphi$,  $\mathscr{A}(\bphi,\cdot)$ and $\mathscr{S}(\bphi,\cdot)$ can be viewed, in a pointwise fashion, as affine linear maps from an affine space $[\E]$  to vector spaces of $H$-valued matrices. Furthermore, these maps have the properties  that if $\Z$ is a section of   $\left(\wh{H}_{m-1} \right)^{d_m}$ then
 \bgE{AJ}  
 \begin{split}
\mathscr{S} (\bphi,  \E+\Z) &= \mathscr{S}(\bphi,  \E) +  \widehat{\J}_{[\bphi]}^{m} ({\Z}) \\
\mathscr{A} (\bphi,  \E+\Z) &= \mathscr{A}(\bphi,  \E) + \J^{m-1,1}([\bphi])(\Z).
\end{split}
\enE
Thus if $M$ is semi-$\J$-nondegenerate, both affine linear maps $\mathscr{S}$ and $\mathscr{A}$ are injective and can be used to identify $[\E]$ pointwise with an affine linear subspace of an appropriate vector space of $H$-valued matrices.

The following result is also trivial.

\bgL{ANat}
Both $\mathscr{S}$ and $\mathscr{A}$ are natural in the sense that
\[  F_*\mathscr{S}(F^* \bphi, \E) =\mathscr{S}(\bphi , F_* \E) \]
for any weak $H$-isometry $F$ with a similar identity for $\mathscr{A}$.
\enL

Before beginning to set up the required algorithm for determining a complement, we establish a technical lemma that we shall need later for vertical rigidity.

\bgL{MPC}
If $M$ is equiregular and semi-$\J$-nondegenerate, then for any $2$-coframe $\bphi$ the map
\[ \text{tr} \widehat{\J}_{\bphi_x}^{2} \colon H_x^{d_2} \to H_x \]
is surjective at all points $x$.
\enL

\pf

Let  $\bphi_x = (\phi^1,\dots,\phi^{d_2})^\top \in \Phi_x^2$. Since $M$ is semi-$\J$-nondegenerate, it follows that $\J^{1,1}[\bphi] \colon H^{d_2} \to M_{d_2 \times d_2} (H)$ is injective and hence
\[ \bigcap\limits_i \ker{ \J^{1,1} ( [\phi^i]_2) } = \{0\}.\]
However, since $ \J^{1,1}([\psi])$ is skew-adjoint for all $\psi \in \wh{V}^2$, we can immediately deduce that
\[ \text{range}(\J^{1,1}( [\phi^1]_2)) +\dots +  \text{range}(\J^{1,1} ( [\phi^{d_2} ]_2)) =H_x .\]
The result then follows from the observation that
\[ \widehat{\J}^{2}_{\bphi_x}(X_1,\dots,X_{d_2}) = 2\sum\limits_{i=1}^{d_2} \J^{1,1}( [\phi^i]_2)(X_i) .\]
\epf

For a weakly non-degenerate equiregular sub-Riemannian manifold this implies that  complements with the required rigidity exist at all points $x \in M$. However, these complements are typically non-unique. To argue that a collection of pointwise complements can be pieced together into a smooth partial complement bundle, we shall need to construct an algorithm for identifying a single pointwise complement from purely sub-Riemannian data and show that the outcome of this algorithm depends smoothly on the parameter $x$. 

The key idea will be to construct partial complements that minimize measures depending on $\mathscr{S}$ and $\mathscr{A}$. As a first step, we establish the following pair of lemmas.

\bgL{minS}
Suppose $M$ is a semi-$\J$-nondegenerate equiregular sub-Riemannian manifold and  $\bphi$ is a smooth local section of  $\Phi^m$ with a dual $(m,m-1)$-frame field $\E$. The function $\|\mathscr{S}(\bphi_x,\cdot)\|$ attains its minimum on $ [\E]$ at a unique $(m,m-2)$-frame field $\F$. This minimizing frame field depends smoothly on $x\in M$.

Furthermore the under the orthonormal change $f \bphi$, the minimal section becomes $\F f^\top$. So this defines a unique, smooth local section of $\Omega^m_{m-2}$.

 \enL

\pf Recall that $[\E]$ is an affine space with underlying vector space $\wh{H}_{m-1}$. The semi-$\J$-nondegeneracy condition then implies that  pointwise $\mathscr{S}$ is an injective, affine linear embedding into $M_{d_m \times d_m}(H)$. Basic linear algebra implies that any affine linear subspace of an inner product space admits a unique element of minimal norm. This minimal element can be constructively found using orthogonal projections and so will depend smoothly on $x \in M$. 

The final part follows easily from \rfE{Sf}.

\epf

A virtually identical argument, using \rfE{Af} in place of \rfE{Sf} proves the following.

\bgL{minA}
Suppose that $M$ is a semi-$\J$-nondegenerate equiregular sub-Riemannian manifold, that $m>k$, that $\bphi$ is a $k$-coframe field and that  $\E$ is a $(m,k-1)$-frame field $\E$ with  $\bphi(\E)=0$. The function $\|\mathscr{A}(\bphi_x,\cdot)\|$ attains its minimum on $ [\E]$ at a unique $(m,k-2)$-frame field $\F$. This minimizing frame field depends smoothly on $x\in M$.

Furthermore the under the orthonormal changes $f_k \bphi$, $\E f_m^\top$ , the minimal section becomes $\F f_m^\top$. So this defines a unique, smooth local section of $\Omega^m_{k-2}$.
\enL

We are now finally in a position to establish our main technical result. 

\bgL{inductive}
Suppose that $0 < m < r$ and local sections $[\ee]^r_{m},[\ee]^{r-1}_{m},\dots,[\ee]^{m+1}_{m}$ of $\Omega^r_m,\dots,\Omega^{m+1}_m$ have been chosen. Then there is an algorithmic method that uniquely determines sections $[\ee]^r_{m-1}, \dots, [\ee]^{m+1}_{m-1},[\ee]^m_{m-1}$ of $\Omega^r_{m-1},\dots,\Omega^{m}_{m-1}$ such that pointwise $[\ee]^j_{m-1} \in [\ee]^{j}_{m}$ for $j=m+1,\dots,r$.
\enL

Before we give a proof, we remark that the algorithm increases the number of sections by one. The additional section is uniquely determined but has no pointwise inclusion condition.

\pf

First we define
\[ V^{m+1} = \left\{ \psi \in H_{m}^o \colon \psi([\ee]^r_{m}) = \dots = \psi( [\ee]^{m+2}_{m}) =0 .\right\}\]
This space is pointwise isomorphic to $\wh{V}^{m+1}$. Choose a representative $\E^{m+1}_m$ of $[\ee]^{m+1}_m$ and a dual coframe $\bphi^{m+1}$ whose entries are contained in $V^{m+1}$. Now we can apply \rfL{minA} to $\bphi^{m+1}$ and the chosen sections one by one. The only exception is $[\ee]^{m+1}_m$ where  we apply \rfL{minS} instead.  This produces the required nested sequence of sections $[\ee]^r_{m-1}, \dots, [\ee]^{m+1}_{m-1}$.

To produce the additional section $[\ee]^m_{m-1}$, we now set
\bgE{Vm}  V^{m} = \left\{ \psi \in H_{m-1}^o \colon \psi([\ee]^r_{m-1}) = \dots = \psi( [\ee]^{m+1}_{m-1}) =0 .\right\}\enE
Choose any $m$-coframe field with entries in $V^m$ and $[\ee]^m_{m-1}$ be the equivalence class in $\Omega^m_{m-1}$ of any dual $(m,m-1)$-frame field.

\epf

We can now prove our main result.\\

\noindent {\em Proof of Theorem A.}  

The proof runs by backwards induction using \rfL{inductive}.  The goal is to determine a unique, intrinsically defined, section $\ee^k_0$ of $\Omega^k_{0}$ for all $2 \leq k \leq r$. 

We begin with the observation that, at the final step, sections of  $\wh{V}^r= (H_{r-1})^o $ are genuine differential forms. Thus all sections of $\Phi^r$ are equivalent up to an orthonormal change.  Choose  a local section $\bphi^{r}$ of $\Phi^r$. This determines a unique local  section $\E^r_{r-1}$ of $\Omega^r_{r-1}(\bphi)$. Under the orthonormal change $\bphi^{r} \mapsto f \bphi^{r}$, we see $\E^r_{r-1} \mapsto \E^r_{r-1} f^\top$ so this determines a unique local section $[\ee]^r_{r-1} \in \Omega^r_{r-1}$.

This establishes a canonical choice of sections as in \rfL{inductive} for the case $m=r-1$.  We now successively apply \rfL{inductive} ending with the case $m=3$. After this application, we have derived canonical sections
\[ [\ee]^r_1,\dots,[\ee]^2_1.\]
For the final step, we partially apply \rfL{inductive} to produce sections
\[ [\ee]^r_0,\dots,[\ee]^3_0.\]
However, we employ a different minimization procedure on $[e]^2_1$ and do not produce the additional section. For the final step of identifying a section $[\ee]^2_0$ pointwise contained in $[e]^2_1$, we first compute a section $R$ of $H$, by
\[ R =  \sum\limits_{m=3}^r \text{tr} \mathscr{S} (\bphi^m,\E^m_0) \]
where $\bphi^m$ is any $m$-coframe field with entries in the bundles $V^m$ defined by \rfE{Vm} and  $\E^m_0$ is the representative of $[\ee]^m_0$ dual to $\bphi^m$. It is easy to check that $R$ is independent under orthonormal changes of each $\bphi^m$.

Next, we let $\bphi^2$ be any $2$-coframe field with entries in $V^2$ and pick $\E^2_1$ to be the representative of $[\ee]^2_1$ dual to $\bphi^2$.  Now, we can define $W \subset \E^2_1$ to be the collection of $(2,0)$-frame fields  $\F$ inside $\E^2_1$ with the property that
\[ \text{tr} \mathscr{S} (\bphi^2,\F) = - R.\]
It follows easily from \rfL{MPC}  that $W$ is non-empty and has  constant rank as an affine subbundle of $\E^2_1$. We then set $\E^2_0$ to be the section of $W$ minimizing $\| \mathscr{S}(\bphi^2,\cdot) \|$. It is easy to check that the equivalence class $[\ee]^2_0$ of $\E^2_0$ is independent of an orthonormal change in $\bphi^2$. 

We now locally define each partial complement $V_m$ to be the span of any representative of $[\ee]^m_0$. The uniqueness of minimizers at each step ensures that these local complements piece together smoothly.  Furthermore it is clear that if we form an orthonormal frame by concatenating representative frames from each $[\ee]^m_0$ the condition of \rfD{srigid} is satisfied. Hence the complement produced form this algorithm is $V$-rigid.
 
To show that the partial complements so produced are all natural, we first remark that all the local bundles used are natural as are the measures being minimized. Naturality follows easily from the fact that the semi-$\J$-nondegenerate condition implies that all the minimizing sections are unique.

\begin{flushright}
$\square$
\end{flushright}

We recall that the complement defined by this algorithm is known as the {\em minimal rigid complement}. The name derives from the minimization procedures that form the core of the algorithm defining the complement.

To illustrate the algorithm, we run it on a simple but non-trivial example.

\bgX{Engelish}
Consider a manifold $M$ with tangent bundle spanned by global, linearly independent vector fields $X_1,X_2,T,S_1,S_2$ where the only non-zero Lie bracket relations are
\[ [X_1,X_2]=T, \; [X_1,T]=S_1, \; [X_2,T]=S_2, \; [X_1,S_1]=3S_2. \]
 Let $\phi^1,\phi^2,\tau,\sigma^1,\sigma^2$ denote the dual frame. If $X_1,X_2$ are declared to be an orthonormal frame for $H$, then it is easy to see that $M$ is equiregular, step $3$ and that the collections $\{T\}$  and $\{S_1,S_2\}$ represent equivalence classes that provide orthonormal frames for $\wh{H}_2$ and $\wh{H}_3$ respectively. Since the bracket structures are identical at every point of $M$, we can do the computations globally.

Set ${\boldsymbol \sigma} = (\sigma^1,\sigma^2)^\top$ then 
\begin{align*} \widehat{J}^{2,1}_{\boldsymbol \sigma} (a[T],b[T]) &= \bp 2aX_1 & bX_1+aX_2 \\ aX_2+bX_1 & 2bX_2 \ep \\
 \widehat{J}^{1,1}_{\tau} (aX_1+bX_2) &=\bp  -2aX_2 +2bX_1\ep. \end{align*}
 Both maps are injective so $M$ is semi-$\J$-nondegenerate.
 
 Now let $\E^3_2 = ([S_1],[S_2])$ be a $(3,2)$-frame dual to ${\boldsymbol \sigma}$. This represents the initial choice in the algorithm. To apply the first iteration of \rfL{inductive}, we compute
 \[ \mathscr{S}({\boldsymbol \sigma},(S_1+ a[T],S_2+ b[T]) )=  \bp 2aX_1 & (b+3)X_1+aX_2 \\ aX_2+(b+3)X_1 & 2bX_2 \ep .\]
 Thus 
 \[ \|\mathscr{S} ({\boldsymbol \sigma},(S_1+ a[T],S_2+ b[T]) )\|^2 = 6a^2+6b^2+12b+18.\] This is minimized when $a=0$, $b=-1 $, so the minimal $(3,1)$-frame $\E^3_1=([S_1],[S_2-T])$. The bundle $V^2$ is then spanned by unit length form $\tau+\sigma^2$ and the additional $(2,1)$- frame is just $\E^2_1 =[T]$. 
 
 For the second iteration of \rfL{inductive}, we compute
 \[ \mathscr{A}(\tau+\sigma^2, S_1+aX_1+bX_2, S_2 -T + cX_1+dX_2) = \bp -aX_2 +(3+b)X_1\\ -cX_2+dX_1 \ep, \]
 which has minimal norm when $a=c=d=0$ and $b=-3$. 
 
 Now 
 \[ R= \text{tr} \mathscr{S}({\boldsymbol \sigma}, (S_1,S_2-T)) = -2X_2 \]
 and since
 \bgE{XE}  \mathscr{S}(\tau+\sigma^2,T+a  X_1 + b X_2) = -2a X_2 + 2b X_1 \enE
we see that $W = T-X_1$ as a $0$-rank affine subbundle of $[T]$.  Since $W$ is $0$-rank, $T-X_1$ must be required minimal section.

Thus the minimal rigid complement has $V_2$ is spanned by $T-X_1$ and $V_3$ spanned by $S_1-3X_2,S_2-T$.
\enX


We conclude this section with the remark that the minimal rigid complement need not be the only natural complement. Indeed, if in the proof of Theorem A we did not restrict our attention to the affine subspace $W$ at the last step, but instead used a final full application of \rfL{inductive}, we would produce an alternative natural complement. In Example \ref{X:none:Engelish} this would mean that we are minimizing in \rfE{XE} with no restrictions on $a$ and $b$, This would change $V_2$ to the span of $T$ instead. 
As we shall see in the next section, the vertical rigidity condition is very useful in applications and so we have chosen an algorithm that produces such a complement. However, it is very possible that this alternate complement could be more useful in other applications. 

With stronger conditions on $M$, it is also likely that natural complements with more refined behavior could be constructed. We do however make the remark that the semi-$\J$-nondegenerate condition is necessary for uniqueness within any such minimization procedure and so will be required for any algorithmic constructions based on the methods outlined here.

\section{Observations and applications}\label{S:obs}

\subsection{Sub-Riemannian Isometries} In \cite{Hladky6}, the author gave an in depth study of the isometries of a complemented sub-Riemannian manifold. Weak $H$-isometries that preserved the complement were called (strong) $H$-isometries.  Now Theorem A provides a semi-$\J$-nondegenerate, equiregular sub-Riemannian manifold with a natural complement. Thus for these manifolds, all weak $H$-isometries must preserve this complement and so the groups $\text{Iso}^*(M)$ and $\text{Iso}(M)$, of smooth weak and strong $H$-isometries respectively,  are identical. Thus the results of  \cite{Hladky6} which held for $\text{Iso}(M)$ must also hold for $\text{Iso}^*(M)$ under these conditions. In particular, from Theorem 3.8 and Corollary 3.9 in \cite{Hladky6}, we can now immediately obtain Theorem B where sub-Riemannian isometries are interpreted as weak $H$-isometries. The results of \cite{Hladky6} can also be used to greatly improve the dimension bound for particular examples. The improvements are difficult to describe for general case and depend closely on the bracket structure. Hence we won't discuss them here.

\subsection{Connections and sub-elliptic operators} Any global rather than local tensorial analysis on a sub-Riemannian manifold requires defining some notion of covariant differentiation. This in turn requires a choice of complement. Usually these covariant derivatives have been defined only on specialized and restrictive categories of sub-Riemannian manifolds. However, in \cite{Hladky4}, a canonical connection was introduced for a sub-Riemannian manifold with a choice of complement and Riemannian extension. In that paper, the Riemannian extension could be chosen arbitrarily provided $H$ and $V$ were orthogonal, but the covariant derivative of a horizontal vector field was invariant under this choice.  However, from our earlier observations, any complement to an equiregular sub-Riemannian manifold carries a canonical Riemannian extension. An intrinsic choice of complement, as offered by Theorem A, now means that, for semi-$\J$-nondegenerate, equiregular sub-Riemannian manifolds, we have a canonical connection defined on the whole tangent bundle.

This connection can essentially be defined as follows: 
\begin{itemize}
\item For $A,B $ sections of the same partial complement $V_m$, $\nabla_A B$ is the projection onto $V_m$ of the Levi-Civita covariant derivative with respect to the canonical metric extension.
\item For  a section $A$ of $V_m$ and $k \ne m$, $\nabla_A$ essentially acts on sections of $V_k$ as the skew-symmetric part of the operator $B \mapsto [A,B]_{(k)}$, where again the subscript $(k)$ denotes projection onto $V_k$. More precisely, if $U_1,\dots,U_m$ is an orthonormal frame for $V_k$, then  $\nabla_A U_i$ is defined by 
\[ \aip{\nabla_A U_i}{U_j}{} =\frac{1}{2}  \left( \aip{[A,U_i]}{U_j}{} - \aip{[A,U_j]}{U_i}{}\right).\]
\end{itemize}
The reader is referred to \cite{Hladky4} for more details. It should be noted that the torsion $T$ for this connection is decidedly non-zero. However for $A,B$ sections of $V_m$, $T(A,B) $ is orthogonal to $V_m$ and for $B \in V_k$ with $k\ne m$ the operator $T(B,\cdot)_{(m)}  \colon V_m \to V_m$ is symmetric.

Various conditions on the torsion have been important in the literature (see for example \cite{BaudoinGarofalo}, \cite{Hladky6},  \cite{Hladky5} ,\cite{Hladky4} ). Here we shall focus on two, $V$-normal and $V$-rigid.

\bgD{VRVN}
The connection above is $V$-normal if $T(H,V_m)$ is orthogonal to $V_m$ for each $m$.  The connection is $V$-rigid if whenever $\{U_i\}$ is an orthonormal frame for $V$,
\[ \sum\limits_i \aip{T(H,U_i)}{U_i}{} =0 .\]
\enD

Since the connection is uniquely defined by the complement, we shall refer to a complement as being $V$-rigid or $V$-normal. It is clear from Theorem A that for a semi-$\J$-nondegenerate, equiregular manifold the determined complement is always $V$-rigid. Furthermore, an easy  examination of semi-$\J$-nondegenerate condition and the proof of Theorem A shows that if a $V$-normal complement exists, it is  unique, natural, and must be the complement constructed in Theorem A.

The importance of $V$-rigid complements comes from two main sources. First, given a complement and the associated connection, we can define the horizontal Laplacian by
\bgE{Lap1} \triangle_\tH \tau = \nabla_{E_i } \nabla_{E_i} \tau -  \nabla_{\nabla_{E_i} E_i} \tau \enE
where $\{E_i\}$ is any local orthonormal frame for $H$. In the presence of a canonical complement and Riemannian metric extension, there is then a canonical sub-elliptic operator associated to $M$ in the same way that the Laplace-Beltrami operator is associated to a Riemannian manifold. A natural question is to ask what properties of the Laplace-Beltrami operator these horizontal Laplacians share. First, we note the existence of a natural volume form.

\bgL{VF}
An orientable equiregular sub-Riemannian manifold of any step size $r$ has a canonical, natural volume form.
\enL

\pf
For $1 \leq m \leq r$, let $\{ \psi^1_m,\dots,\psi^{d_m}_m\}$ be a collection of local $1$-forms in $H_{m-1}^o$ such their representatives in $\wh{V}^m$ form a local orthonormal frame.  We note then note that
\[ \Psi_m = \psi^1_m \wedge \dots \wedge \psi^{d_m}_m \]
is uniquely defined modulo $H_m^o$,  up to sign. Thus it is easy to see that
\[ \Psi_1 \wedge \dots \Psi_r \] 
is uniquely defined up to sign. If $M$ is orientable, the sign can then be chosen consistently to yield a canonical volume form. That this volume form is natural is trivial and left to the reader.  

\epf

This volume form is sometimes known as the Popp volume form, see \cite{Barilari}. It is clear from the proof, that this volume for does not depend in any way on the complement. However, a fixed choice of partial complements means that the $1$-forms can be chosen to be explicit orthonormal coframes on each level of the complement.

It can easily be seen (see \cite{Hladky4}) that if $M$ is compact and the complement is $V$-rigid then on functions 
\bgE{Lap2}  \triangle_\tH =  - \nabla^* _\tH \nabla_\tH  \enE
where the $L^2$-adjoint is taken with respect the canonical volume form of \rfL{VF}. Hence in the $V$-rigid case, the horizontal Laplacian is self-adjoint.  

\subsection{Minimal and CMC surfaces} The second area where $V$-rigid complements are important is the study of minimal and constant mean curvature hypersurfaces in sub-Riemannian geometry. From the work of \cite{HP2} and \cite{HP4}, in the $V$-rigid case the first and second variation of perimeter measure can be handled in a relatively straightforward manner. The above hypersurfaces can be characterized by a mean curvature tensor that is closely related to the classical tensor from Riemannian geometry. Without the $V$-rigid condition an extra term appears in all the the characterizing equations and everything becomes much harder to work with (see for example Section 5 of \cite{HP4}). For example, there are no second variation formulas on general categories of non $V$-rigid sub-Riemannian manifolds.  

\section{The step $2$ case}\label{S:two}

The semi-$\J$-nondegenerate condition is least restrictive and easiest to analyze in the step $2$ case. This is largely due to the fact that we need only study a single operator $\J =\J^{1,1} \colon H^o_1 \to \text{Hom}(H,H)$ and so all of the spaces involved consist of genuine vectors and covectors without the need for quotient spaces and equivalence classes. Here also we can define a stronger nondegeneracy condition.

\bgD{JND}
A step $2$ sub-Riemannian manifold is $\J$-nondegenerate at $x$ if for all $\phi \in V_x^* \backslash \{0\}$, the map $\J(\phi)$ is an isomorphism of $H_x$.
\enD

Given an orthonormal frame $\phi^1,\dots,\phi^{d_2}$ for $H_1^o$, looking for a $V$-normal complement amounts to looking for a dual frame $U_1,\dots,U_{d_2}$ such that for all $i,j$
\[ d\phi^i(U_j,X) + d\phi^j(U_i,X) = 0 \]
for all $X \in H$. The vector fields $U_i$ can be interpreted as a natural generalization of the notion of a Reeb vector field for a contact form. In the notation of this paper this is looking for a dual frame $\boldsymbol U$ such that $\mathscr{S} (\bphi,\boldsymbol U) = 0$. This idea appears naturally in the study of quaternionic contact manifolds conducted by Biquard \cite{Biquard} and Duchemin \cite{Duchemin}. Biquard established the existence of a $V$-normal complement for quaternionic manifolds of dimension $4n+3$ with $n>1$. In general, even assuming $\J$-nondegeneracy,  it  is not always possible to construct $V$-normal complements. For quaternionic contact manfolds, this occurs in dimension $7$. Duchemin extension of Biquard's work to dimension $7$ assumes the existence of a $V$-normal complement. The canonical $V$-rigid complements constructed in this paper agree with these $V$-normal complements when they exist and should provide a viable alternative when they do not.

As the semi-$\J$-degenerate condition is somewhat unwieldy, we conclude with some simpler sufficient conditions.  
\bgL{suff}
The following conditions are independently sufficient for $M$ to be semi-$\J$-nondegenerate at $x$,

\begin{itemize}
\item $M$ is $\J$-nondegenerate at $x$.
\item There exists a $\psi \in (H_1^o)_x$ such that $\J(\psi)$ is an isomorphism of $H$.
\item There exists a orthonormal conframe $\bphi = (\phi^1,\dots,\phi^{d_2})^\top$ such that $\J(\phi^1) + \dots + \J(\phi^{d_2})$ is injective on $\ker{\J(\phi^1)} + \dots + \ker{\J(\phi^{d_2})}$.
\end{itemize}
\enL

\pf For simplicity of notation, if we have chosen a coframe $\bphi$, we shall set $J_i = \J(\phi^i)$. Sufficiency of the first condition is trivial and left to the reader. (It also follows from either of the other two.)

For the second condition, we note that by an orthonormal change of frame and a re-scaling, we can assume that $\psi = \phi^1$ for some coframe $\bphi=  (\phi^1,\dots,\phi^{d_2} )^\top$ and so $J_1$ is an isomorphism.  Then if $\widehat{J}_{\bphi}(\Z) =0 $ we clearly must have $J_1 Z_1=0$ and so $Z_1 =0$. But then $0 = J_k Z_1 + J_1 Z_k = J_1 Z_k$ and so clearly $\Z =0$. 

For the third condition, we shall argue the case $d_2>2$ and leave the trivial cases $d_2 \leq 2$ to the reader.  We first set $N = \ker J_1 + \dots + \ker J_{d_2}$ and $J_\Sigma = J_1 +\dots J_{d_2}$. If $\widehat{J}_{\bphi}(\Z) =0 $, then clearly each $Z_i \in \ker J_i \subseteq N$ and for all $1\leq i,j \leq d_2$, $(J_i+J_j)(Z_i+Z_j) =0$. However it is then easy to see that 
\[ J_\Sigma (Z_1 +\dots + Z_{d_2} ) = \sum\limits_{1 \leq i \leq j \leq d_2 } (J_i+J_j)(Z_i+Z_j)  =0.\]
But then by assumption $Z_1 = - (Z_2 + \dots Z_{d_2})$ and for $1< i<j \leq d_2$,
\begin{align*}  (J_i + J_j)v_1 &= -(J_i+J_j) ( Z_2 + \dots+ Z_{d_2} ) \\
&= - (J_i+J_j)( Z_2 + \dots \widehat{Z_i} + \dots + \widehat{Z_j} + \dots + Z_{d_2}) \end{align*}
where $\hat{}$ denotes omission. Thus, we can sum  to find
\begin{align*}  (d_2-1) J_\Sigma Z_1& = - (d_2-3)  \sum\limits_{  1<i, j\leq d_2} J_i Z_j t = - (d_2-3) (J_\Sigma - J_1) ( Z_2 + \dots Z_{d_2} ) \\
&= (d_2-3)(J_\Sigma -J_1)(Z_1)  = (d_2 -3) J_\Sigma Z_1.
\end{align*}
But this clearly implies that $J_\Sigma Z_1 =0$ and so  $Z_1 =0$. That all other $Z_k=0$ follows similarly.


\epf

As a final comment, we note that this implies that all strictly pseudoconvex pseudohermitian manifolds and all quaternionic contact manifolds (including dimension $7$) are semi-$\J$-nondegenerate. In particular Theorem B applies to both categories.

\bibliographystyle{plain}
\bibliography{/Users/robert/Dropbox/References}
 \end{document}